\documentclass{article}
\usepackage{graphicx} % Required for inserting images
\usepackage{amssymb}
\usepackage{amsmath}
\usepackage{bbold}
\usepackage{algorithm}
\usepackage{algpseudocode}
\usepackage{pgfplots}
\pgfplotsset{compat=1.8}
\usepgfplotslibrary{statistics}
\usepackage[T1]{fontenc}
\usepackage[utf8]{inputenc}

\usepackage[colorlinks=true, allcolors=blue]{hyperref}  %% 
\usepackage{cleveref} % For nicer references

\usepackage{amsthm}
\usepackage{xcolor}
\usepackage{multirow}
\usepackage{authblk}

\let\oldtextit\textit 
\renewcommand\emph[1]{\oldtextit{\color{RoyalBlue}#1}}
\definecolor{RoyalBlue}{cmyk}{1, 0.50, 0, 0}
\makeatletter
\def\thanks#1{\protected@xdef\@thanks{\@thanks
        \protect\footnotetext{#1}}}
\makeatother
\theoremstyle{definition}

\newtheorem{theorem}{Theorem}

\newtheorem{lemma}[theorem]{Lemma}

\newtheorem{remark}[theorem]{Remark}
\newtheorem{example}[theorem]{Example}
\newtheorem{conjecture}[theorem]{Conjecture}
\newtheorem{que}[theorem]{Question}

\title{Asymptotic rank bounds: a numerical census}

\author{Kisun Lee\thanks{
\hspace*{-1.8em}School of Mathematical and Statistical Science, Clemson University, 220 Parkway Drive, Clemson, SC 29634 
\\
\texttt{email:}\href{mailto:kisunl@clemson.edu}{kisunl@clemson.edu}, \texttt{ URL:}\href{https://klee669.github.io}{https://klee669.github.io}}}
\date{}

\begin{document}

\maketitle
\begin{abstract}
We systematically compute improved asymptotic rank bounds for tensors. Using numerical implicitization, we implement the geometric framework of Kaski and Micha{\l}ek across all computationally feasible cases. By detecting the absence of low-degree vanishing polynomials on secant varieties, we obtain new asymptotic rank bounds that improve upon the generic border rank bounds. The results provide numerical data supporting Strassen's asymptotic rank conjecture and clarify the computational barriers posed by current numerical methods.
\end{abstract}
\section{Introduction}

The behavior of tensor powers plays a central role in the study of algorithms and algebraic complexity theory \cite{burgisser2013algebraic}. For a fixed tensor $T \in \mathbb{k}^a \otimes \mathbb{k}^b \otimes \mathbb{k}^c$, the growth rate of the rank of its powers $T^{\otimes q}$ reflects the asymptotic complexity inherent in $T$. This viewpoint is classical: Strassen showed that the exponent of matrix multiplication is determined by the asymptotic rank of the specific $2\times 2\times 2$ tensor \cite{strassen1986asymptotic}. Thus, understanding the long-term behavior of tensor powers is tightly connected to key open problems in algebraic complexity.

A geometric approach, building on \emph{border rank} and \emph{secant varieties}, offers a complementary perspective. One studies spaces of tensors to identify the worst-case behavior \cite{bini1979n2,schonhage1981partial}. However, this approach faces a significant challenge: proving bounds on the asymptotic rank of specific tensors is notoriously difficult. While random tensors are known to have maximal border rank, making high rank tensors the ``hay'' in a haystack, constructing an explicit sequence of tensors with maximal border rank has proven largely unsuccessful (a difficulty termed ``finding hay in a haystack'' by Howard Karloff).

This difficulty highlights a contrast between generic behavior and the asymptotic behavior of structured tensors. In this context, a natural class of examples consists of \emph{concise} tensors, those whose flattenings are injective, which in our setting corresponds to formats satisfying $a\leq b\leq c$ and $c<ab$. Another class is that of \emph{tight} tensors, whose support satisfies an additive conservation law (see \cite[Definition~1.1]{conner2021towards}). For such tensors, the generic border rank is typically larger than the dimensions $a,b$ and $c$. In contrast, the \emph{asymptotic rank conjecture} (see \cite[Conjecture 1.4]{conner2021towards} and \cite[Section 13]{wigderson2022asymptotic}) predicts that their asymptotic behavior is far simpler, collapsing to the lower bound determined solely by the format.

\begin{conjecture}[Strassen's asymptotic rank conjecture]
    For any tight and concise tensor $T \in \mathbb{k}^a \otimes \mathbb{k}^b \otimes
    \mathbb{k}^c$,
    \[
    \lim_{q\to\infty} \text{rank}(T^{\otimes q})^{\frac{1}{q}}
    =
    \max\{a,b,c\}.
    \]
\end{conjecture}

An affirmative answer to this conjecture would imply that for structured tensors, the asymptotic rank matches the trivial lower bound. This would provide a universal upper bound on the complexity of tensor problems (e.g.\ matrix multiplication) without the need to analyze each problem individually. For more context, we refer to the introduction of \cite{kaski2025universal} and references therein.

Despite this appealing conjecture, asymptotic rank remains barely understood except for some supporting results \cite{christandl2019barriers,conner2022rank,strassen1988asymptotic}. Even for small formats, determining the scaling of $\text{rank}(T^{\otimes q})$ is challenging because asymptotic behavior cannot be read off from classical invariants like border rank (see \cite{buczynska2014secant,galkazka2023distinguishing}, for instance). Consequently, establishing sharper upper bounds requires new ideas, and so even incremental progress may be significant.

We perform a systematic numerical computation of asymptotic rank bounds for tensors in $\mathbb{k}^a\otimes\mathbb{k}^b\otimes \mathbb{k}^c$ across all formats whose generic border rank is at most $20$. Through this extensive search, we obtain $5$ cases with strictly improved bounds within the current computational limits; see \Cref{tab:results}. Our approach applies numerical implicitization to implement the framework established in \cite{kaski2025universal}. The results provide empirical evidence for asymptotic behavior previously beyond the reach of theoretical approaches and quantify the computational barriers in higher codimension cases.

\Cref{sec:geometric_framework} outlines the framework from \cite[Section 5.1]{kaski2025universal}. \Cref{sec:numerical_implicitization} details the numerical implicitization. Lastly, \Cref{sec:results} presents the experimental results.

\section{Geometric framework}\label{sec:geometric_framework}
In this section, we elaborate on the geometric framework established in \cite{kaski2025universal}, derive asymptotic rank bounds from secant varieties, and illustrate its application with an example.

Consider a field $\mathbb{k}$ and positive integers $a,b,c$. Let $V = \mathbb{k}^a\otimes\mathbb{k}^b\otimes \mathbb{k}^c$. For a positive integer $q$, we analyze the \emph{Kronecker power map} (or the $q$-th \emph{Veronese embedding}):
\begin{equation*}
    K_{q}: V \rightarrow V^{\otimes q}, \quad T \mapsto T^{\otimes q}.
\end{equation*}
The image of this map is contained in the $\mathfrak{S}_q$-invariant subspace $(V^{\otimes q})^{\mathfrak{S}_q}$ of $V^{\otimes q}$. This geometric setting allows us to derive upper bounds for the rank of $K_q(T)=T^{\otimes q}$ for any tensor $T\in V$.

We study the rank of $T^{\otimes q}$ from a linear algebra perspective within this $\mathfrak{S}_q$-invariant subspace $(V^{\otimes q})^{\mathfrak{S}_q}$. This subspace admits a basis $\{T^{(g)}\}$ indexed by integer compositions $g: \Delta \to \mathbb{N}_{\ge 0}$, where $\Delta=[a]\times [b]\times [c]$ and $\sum_{\delta \in \Delta} g(\delta) = q$. Note that $g$ records how many times each index $\delta$ appears in the $q$ tensor factors. The size of this basis is given by ${abc+q-1\choose q}$. With respect to this basis, the map decomposes explicitly as:
\begin{equation}\label{eq:decomposition}
    K_q(T) = T^{\otimes q} = \sum_{g} T^g\, T^{(g)}.
\end{equation}
where the coefficient of $T^{(g)}$ is given by the monomial $T^g:=\prod_{\delta \in \Delta} T_{\delta}^{g(\delta)}$ of degree $q$. %This decomposition reveals that determining the asymptotic rank reduces to analyzing the rank of these explicit basis tensors $T^{(g)}$.
Moreover, by \cite[Proposition 15]{kaski2025universal}, we have 
$\text{span}\{T^{(g)}\} = \text{span}(K_q(V))$. Combined with the subadditivity of tensor rank, this yields the following bounds:
\begin{lemma}\cite[Proposition 15 and Lemma 17]{kaski2025universal} \label{lem:rank_bounds}
We have
\begin{multline*}
{abc+q-1 \choose q} \cdot \max_{g} \text{rank}(T^{(g)})
\geq \max_{T \in (V^{\otimes q})^{\mathfrak{S}_q}} \text{rank}(T)
\\\geq \max_{T \in V} \text{rank}(T^{\otimes q})
\geq \frac{1}{{abc+q-1 \choose q}} \max_{g} \text{rank}(T^{(g)}).
\end{multline*}
\end{lemma}

Recall that the asymptotic rank is defined through the $q$-th roots of the ranks of tensor powers.  Hence, Lemma~\ref{lem:rank_bounds} suggests a strategy for obtaining asymptotic rank bounds: it expresses $\text{rank}(T^{\otimes q})$ using the decomposition \eqref{eq:decomposition}. This decomposition shows that an asymptotic bound is governed by how the ``contributions'' of the terms $T^g T^{(g)}$ behave as $q$ grows. Thus, the problem reduces to determining how many coefficients $T^g$ can be nonzero simultaneously. This can be achieved through the secant variety.

Before turning to the secant variety viewpoint, we recall the following general principle from \cite{kaski2025universal}, which formalizes how non-vanishing of degree $q$ polynomials on a subset yields asymptotic rank bounds.
\begin{theorem}\cite[Theorem 44]{kaski2025universal}\label{thm:core}
    Let $L$ be a fixed subspace in $V=\mathbb{k}^a\otimes \mathbb{k}^b\otimes\mathbb{k}^c$, and $Y\subset L$ be a subset with the property that for all $T\in Y$, the asymptotic rank is at most $r$. Suppose that there is no homogeneous polynomial on $L$ of degree $q$ that vanishes on $Y$. Then, every tensor in $L$ has asymptotic rank at most
    \[r{\dim L+q-1\choose q}^{\frac{1}{q}}.\]
\end{theorem}
We present one application of this principle.

\begin{example}\cite[Example 45 revisited]{kaski2025universal}
We consider the space $V=\mathbb{C}^7\otimes\mathbb{C}^7\otimes\mathbb{C}^7$.
    \begin{enumerate}
        \item \textbf{Identify the known upper rank:} The generic border rank in $V$ is $19$. This provides an upper bound of $19$ for the asymptotic rank. We aim for a bound lower than $19$.
        
        \item \textbf{Geometric setup:} We focus on the secant variety $\sigma_{18}(V)$. (We write $\sigma_r(V)$ for the $r$-th secant variety of the variety of rank-one tensors in $V$.) Note that it forms a hypersurface of degree at least $187000$ \cite{hauenstein2013equations}.

        \item  \textbf{Subspace restriction:} Consider a generic line $L$. Since $\sigma_{18}(V)$ is a hypersurface, we have $Y=L\cap \sigma_{18}(V)$ with $|Y|\geq 187000$.

        \item \textbf{Polynomial absence:} We consider homogeneous polynomials on the line $L$ of degree $186999$ (one less than the degree bound). The dimension of the space of such polynomials (which is spanned by monomials) is $187000$. Since a non-zero univariate polynomial of degree $186999$ can have at most $186999$ roots, and the intersection contains at least $187000$ points, no non-zero polynomial of degree $186999$ vanishes on $Y$. This non-vanishing condition implies that the linear span of $K_{186999}(\sigma_{18}(V))$ is the invariant subspace $(V^{\otimes 186999})^{\mathfrak{S}_{186999}}$. Hence, the rank of $T^{\otimes 186999}$ for any $T\in V$ is bounded above by the maximal rank among $S^{\otimes 186999}$ with $S\in\sigma_{18}(V)$.
        
        \item \textbf{Bound computation:} For any $S\in \sigma_{18}(V)$, note that 
        \[\text{rank}(S^{\otimes 186999})\leq 18^{186999}.\]
        Since $(V^{\otimes 186999})^{\mathfrak{S}_{186999}}$ has $187000$ basis elements, the asymptotic rank is bounded by
        \[
    18 \cdot 187000^{\frac{1}{186999}} < 18.001169.
        \]
    This yields an improved upper bound on the asymptotic rank in $V$.
    \end{enumerate}
    This example illustrates the use of \Cref{thm:core}, and provides a bound that is sharper than the one obtained in \cite[Example 45]{kaski2025universal}. This example serves as a prototype for our numerical experiments, in which we approximate the same non-vanishing condition by numerically sampling the secant variety.
\end{example}

Our approach relies on determining both the dimension and degree of the secant variety, and on proving the non-existence of nonzero polynomials of degree at most $q$ vanishing on the secant variety. Note that this can be (computationally) done in practice by sampling points from the intersection $L\cap \sigma_{r}(V)$. We carry this out using numerical methods.

\section{Numerical implicitization of secant varieties}\label{sec:numerical_implicitization}

To describe how the geometric framework from the previous section is implemented numerically, we give an overview of the \emph{numerical implicitization}. It is used to compute the intersection $L\cap \sigma_r(V)$. Our goal is not necessarily to compute the exact degree of a secant variety, but to sample sufficiently many points in $L\cap \sigma_r(V)$ to guarantee the absence of low-degree vanishing polynomials needed for the asymptotic rank bounds. From now on, we assume that $\mathbb{k}=\mathbb{C}$ as we conduct numerical computation.

A classical \emph{witness set} \cite[Section 13.3]{sommese2005numerical} for an irreducible variety $X$ consists of a polynomial system defining $X$, together with its intersection with a generic linear slice $L$ of complementary dimension. It encodes both the dimension and degree of $X$ using a finite set of points. However, if the defining equations are inaccessible, it is hard to construct classical witness sets.

Instead, when $X$ is given through a parametrization as a projection of a polynomial system, we use a \emph{pseudowitness set} \cite{hauenstein2010witness,hauenstein2013membership}. It plays the same role at the level of the image, enabling the sampling of points on $X$ and the computation of its dimension and degree, without defining equations. This approach is called the \emph{numerical implicitization}.

Recall that $V=\mathbb{C}^a\otimes \mathbb{C}^b\otimes \mathbb{C}^c$.
For the $r$-th secant variety $\sigma_r(V)$, we define the parametrization map on the affine space:
$$\begin{array}{rccc}
\Phi_r : & \left(\mathbb{C}^{a}\times\mathbb{C}^{b}\times\mathbb{C}^{c}\right)^{ r} &\rightarrow& V \\
& ((a_1, b_1, c_1), \dots, (a_r, b_r, c_r)) &\mapsto & \sum\limits_{i=1}^{r} a_i\otimes b_i\otimes c_i
\end{array}$$
Rather than treating $\sigma_r(V)$ as an implicit algebraic set defined by equations, we work with this explicit parametrization. 

We work in the domain of $\Phi_r$ with variables $(a_1,\dots a_r;b_1,\dots,b_r;c_1,\dots, c_r)$. To determine the dimension of $\sigma_r(V)$, we choose $\ell$ generic linear functionals $\lambda_1,\dots,\lambda_\ell\in V^\star$ and impose the pulled-back equations $\lambda_j\circ\Phi_r(u)=0$ on the parameter variables. The number $\ell$ is increased until the projection onto $V$ becomes zero-dimensional,
so that $L\cap\sigma_r(V)$ is finite. In this case, the intersection $L\cap\sigma_r(V)$ is parametrized by the solutions of a polynomial system in the variables $(a_1,\dots a_r;b_1,\dots,b_r;c_1,\dots, c_r)$.

Once a single solution is known, the remaining intersection points are recovered by homotopy continuation \cite[Chapter 2]{sommese2005numerical}. We deform the linear constraints $\lambda_j$ along a loop in the space of linear slices and track the corresponding solution path. When the constraints
return to their initial values, the endpoint typically lands at a different solution.
This permutation of solutions is the \emph{monodromy action}
\cite{duff2019solving,sommese2001using}. By iterating such loops, we generate new solutions. When all solutions to the system are found, the number of unique points in their images equals the degree of $\sigma_r(V)$. General procedures of numerical implicitization are implemented in \cite{chen2019numerical}.

\begin{remark}
In the actual implementation of the above procedure, we do not expand the ambient tensor coordinates, but instead preserve the multilinear structure of the rank $r$ parametrization $\Phi_r$ for the efficiency of homotopy continuation. Each rank-one summand is written in affine charts as
\[
(a_i,1) \otimes (b_i,1) \otimes c_i,
\quad
a_i\in\mathbb{C}^{a-1},\; b_i\in\mathbb{C}^{b-1},\; c_i\in\mathbb{C}^{c},
\]
so that the parametrized tensor is
\[
T(a,b,c) = \sum_{i=1}^r (a_i,1) \otimes (b_i,1) \otimes c_i \in V.
\]

A generic linear subspace $L\subset V$ of dimension $\ell$ (where $\ell$ is the codimension of $\sigma_r(V)$) is represented in as $L = \{A t + B \mid t\in\mathbb{C}^{\ell}\}$. The image-space slicing equations are imposed by solving \begin{equation}\label{eqref:system}
    T(a,b,c) - (A t + B)=0
\end{equation} for the unknowns $(a_1,\dots, a_r;b_1,\dots, b_r;c_1,\dots, c_r)$ and $t=(t_1,\dots, t_\ell)$. We initialize the parameters $A\in\mathbb{C}^{abc\times \ell}$ and $B\in \mathbb{C}^{abc}$ generically and vary them along monodromy loops. This produces the pulled-back nonlinear equations $\lambda_j\circ\Phi_r(a,b,c)=0$ without forming the coordinates of $T(a,b,c)$ as independent variables.

To eliminate the remaining degrees of freedom in the parameter space, additional generic linear equations may be imposed directly on the unknowns $u=(a_1,\dots, a_r;b_1,\dots, b_r;c_1,\dots, c_r)$, i.e.\ \[
H u - H u_0 =0,
\quad H\in\mathbb{C}^{m\times r((a-1)+(b-1)+c)} \quad\text{and}\quad u_0 \in \mathbb{C}^{r((a-1)+(b-1)+c)},
\]
where $m$ equals the generic fiber dimension of $\Phi_r$. These \textit{fiber slices} do not further restrict the image $L\cap\sigma_r(V)$ but serve only to obtain a square system suitable for homotopy continuation. 
\end{remark}

This procedure provides a finite subset of $L\cap\sigma_r(V)$. We sample a large enough subset to prove the non-existence of vanishing polynomials on $L\cap\sigma_r(V)$.

\section{Experimental results}\label{sec:results}

Based on the discussions above, we compute asymptotic rank bounds using \Cref{thm:core}. For each triple $(a,b,c)$ satisfying the conciseness condition ($a \le b \le c$ and $c < ab$), we compute the \emph{generic border rank} of $V=\mathbb{C}^a \otimes \mathbb{C}^b \otimes \mathbb{C}^c$ using numerical implicitization. The generic border rank is the minimal integer for which the associated secant variety fills the ambient space. We focus on the cases where this generic border rank is strictly larger than $\max\{a,b,c\}$.

For integers $r$ strictly below the generic border rank, we consider the secant variety $\sigma_r(V)$ and compute its dimension and degree using the same numerical procedure. As the systems grow rapidly with $(a,b,c)$, we restrict to cases that remain computationally feasible. In particular, we examine all concise formats with $r<20$ (that is, the case with the generic border rank at most $20$).

All experiments were conducted using \texttt{HomotopyContinuation.jl}~\cite{breiding2018homotopycontinuation} on a high-performance computing cluster. Computations were carried out on a Genoa-generation AMD CPU node, with each job allocated 16 CPU cores and 64 GB of RAM, and a maximum wall time of 72 hours. 

\subsection{Codimension \texorpdfstring{$1$}{1} cases}

The results of codimension $1$ cases are summarized in \Cref{tab:results}. With more computational resources, one may obtain sharper bounds by sampling more points except for the case $\sigma_8(\mathbb{C}^3\otimes \mathbb{C}^5\otimes \mathbb{C}^7)$. For the cases of the form $\sigma_{3n+1}(\mathbb{C}^{3}\otimes\mathbb{C}^{2n+1}\otimes \mathbb{C}^{2n+1})$, all of them have degree at most $6n+3$ (see~\cite[Section~10.3]{michalek2021invitation}), which do not induce improved asymptotic rank bounds. The next promising case is $\sigma_{22}(\mathbb{C}^{5}\otimes\mathbb{C}^{7}\otimes \mathbb{C}^{17})$, but the implementation failed to track a monodromy loop from the initial solution.  %All \texttt{Julia} scripts are available at:

\begin{table}[ht]
\centering
\renewcommand{\arraystretch}{1.20}
\noindent\makebox[\textwidth]{\begin{tabular}{c||c|c|c|c}  
$(r,a,b,c)$  & \# vars & \# params &  degree & \begin{tabular}{@{}c@{}} asymptotic\\ rank bound\end{tabular} \\
\hline
 $(3n+1,3,2n+1,2n+1)$  & $3(2n+1)^2$ & $6(2n+1)^2$ & 
$6n+3$ & 
N/A \\
 $(8,3,5,7)$  & $105$ & $210$ & 
$105$ \cite{hauenstein2013equations} & 
$<8.366128$ \\
$(17,4,7,14)$  & $392$ & $784$ & $\geq 1229$
 
 & $<17.098769$ \\
$(17,6,6,9)$  & $324$ & $648$ & $\geq 3601$
 
 & $<17.038715$ \\
$(18,7,7,7)$ & $343$ & $686$ & $\geq187000$\cite{hauenstein2013equations} 
 & $<18.001169$ \\
$(19,5,8,10)$ & $400$ & $800$ &  
$\geq 3638$ & $<19.042882$
\end{tabular}}
\smallskip
\caption{Summary of numerical results for concise tensor formats $(a,b,c)$ with $\sigma_{r}(V)$ of codimension $1$. The number of variables and parameters from the system \eqref{eqref:system} are listed in the \# vars and \# params columns. The degree column lists lower bounds of $\deg(\sigma_r(V))$ obtained via numerical implicitization. N/A indicates that the degree is insufficient to improve upon the generic bound.}\label{tab:results}
\end{table}

\subsection{Codimension \texorpdfstring{$2$}{2} and \texorpdfstring{$3$}{3} cases}

When the codimension exceeds $1$, the size of $L \cap \sigma_r(V)$ is no longer sufficient to prove the absence of a vanishing polynomial, unless the arithmetic Cohen-Macaulay condition is imposed. In this case, establishing a non-vanishing condition requires \emph{polynomial interpolation}: sample sufficiently many points and test if all degree $q$ polynomials vanish on the intersection. This requires verifying whether the interpolation matrix has full rank.

For example, the smallest codimension $2$ case is $\sigma_9(\mathbb{C}^4\otimes\mathbb{C}^4\otimes\mathbb{C}^8)$. 
%(smaller codimension $2$ cases than this do not provide enough points for the interpolation). 
To improve the asymptotic rank bound using \Cref{thm:core}, one would need to show the non-existence of a vanishing polynomial of degree at least $q=76$. The corresponding interpolation matrix is indexed by all homogeneous monomials in three variables of degree $76$, and has size at least $3003 \times 3003$. Such matrices are numerically ill-conditioned and require preconditioning techniques, as in \cite{griffin2014numerical}, after high-precision refinement of the sampled points. Despite these, numerical instability prevented reliable rank verification. We report the lower bounds on $\deg(\sigma_r(V))$ via numerical implicitization instead in \Cref{tab:results-codim23}.

\begin{table}[ht]
\centering
\renewcommand{\arraystretch}{1.20}
\noindent\makebox[\textwidth]{
\begin{minipage}{0.64\textwidth}
\centering
\begin{tabular}{c||c|c|c} 
$(r,a,b,c)$  & codim &  degree & minimal $q$ \\
\hline
$(9,4,4,8)$  & $2$ & $\geq 30005$ & $76$ \\
$(10,3,6,9)$  & $2$ & $\geq 78589$ & $87$ \\
$(11,3,7,9)$  & $2$ & $\geq 23724$ & $98$  \\
$(13,5,6,7)$  & $2$ & $\geq 3105$ & $121$ \\
$(14,5,6,8)$  & $2$ & $\geq 1767$ & $132$ \\
$(18,4,8,13)$  & $2$ & $\geq 1057$ & $180$ \\
$(19,5,7,12)$  & $2$ & $\geq 2333$ & $192$ \\
\end{tabular}
\end{minipage}
\hfill
\begin{minipage}{0.64\textwidth}
\centering
\begin{tabular}{c||c|c|c} 
$(r,a,b,c)$  & codim &  degree & minimal $q$ \\
\hline
$(7,4,4,5)$  &$3$ & $44000$ \cite{hauenstein2013equations}  & $88$\\
$(9,4,5,6)$ &$3$ & $\geq 33634$ & $120$ \\
$(11,4,6,7)$  &$3$ & $\geq 8625$ & $154$\\
$(12,3,7,11)$ &$3$ & $\geq 2888$ & $171$\\
$(13,4,7,8)$  &$3$ & $\geq 2503$ & $189$\\
$(14,3,9,11)$  &$3$ & $\geq 879$ & $207$\\
$(15,4,8,9)$  &$3$ & $\geq 842$ & $225$\\
$(17,4,9,10)$  &$3$ & $\geq 327$ & $262$\\
$(17,5,6,12)$  &$3$ & $\geq 317$ &$262$ \\
$(19,4,10,11)$  &$3$ & N/A & $299$ \\
\end{tabular}
\end{minipage}
}
\smallskip
\caption{Numerical data for concise tensor formats $(a,b,c)$ for which the secant variety $\sigma_r(V)$ has codimension $2$ or $3$. The parameter $q$ is the minimal degree required for the improved asymptotic rank bound.  N/A indicates the failure to track a monodromy loop from the initial solution.}
\label{tab:results-codim23}
\end{table}

When the codimension exceeds $3$, the size of the interpolation matrix becomes prohibitively large. For instance, the smallest codimension $4$ case is $\sigma_5(\mathbb{C}^3 \otimes \mathbb{C}^4 \otimes \mathbb{C}^4)$, for which the minimal degree required for applying \Cref{thm:core} is $q = 80$. The corresponding interpolation matrix has dimension exceeding $1.9\times 10^{6}$, which is infeasible in practice.

\section{Concluding remark}

From the experimental results, we observe that all cases are non-defective except for the family $\sigma_{3n+1}(\mathbb{C}^{3}\otimes\mathbb{C}^{2n+1}\otimes \mathbb{C}^{2n+1})$. In particular, for codimension 1 cases, the degree of the secant variety is consistently large enough to improve the geometric bounds. This leads us to the following question: 
\begin{que}
Does every non-defective case admit a strict improvement via the framework of \Cref{thm:core}?    
\end{que}
The challenge lies in the trade-off between the rank reduction and the geometric penalty in \Cref{thm:core}. As $r$ decreases, the codimension of $\sigma_r(V)$ increases, inflating the binomial coefficient ${\dim L +q - 1\choose q}$. For an improved bound to exist, the degree of $\sigma_r(V)$ (and thus the admissible $q$) must grow sufficiently fast to offset this penalty. Our exhaustive search, conducted regardless of codimension, found no counterexamples: in every non-defective case with $r\leq10$, the degree was large enough to yield an improved bound. While this suggests a positive answer to the question, since the degree alone does not directly determine $q$ in higher codimensions, verifying this requires improvements in numerical techniques.

\section*{Acknowledgement}
We thank the organizers of the AGSTA Summer School on Tensors 2025, where this project began. We are grateful to Mateusz Micha{\l}ek for introducing the problem and encouraging the project, and for helpful discussions. We also thank Maciej Ga{\l}{\k{a}}zka for insightful discussions. Finally, we appreciate the help of Paul Breiding, Timothy Duff, and Jonathan Hauenstein with numerical computations.

\bibliographystyle{abbrv}
\bibliography{ref}

\begin{thebibliography}{10}

\bibitem{bini1979n2}
D.~Bini, M.~Capovani, F.~Romani, and G.~Lotti.
\newblock ${O}(n^{2.7799})$ complexity for $n\times n$ approximate matrix
  multiplication.
\newblock {\em Information processing letters}, 8(5):234--235, 1979.

\bibitem{breiding2018homotopycontinuation}
P.~Breiding and S.~Timme.
\newblock Homotopy{C}ontinuation.jl: A package for homotopy continuation in
  {J}ulia.
\newblock In {\em International congress on mathematical software}, pages
  458--465. Springer, 2018.

\bibitem{buczynska2014secant}
W.~Buczy{\'n}ska and J.~Buczy{\'n}ski.
\newblock Secant varieties to high degree {V}eronese reembeddings,
  catalecticant matrices and smoothable {G}orenstein schemes.
\newblock {\em Journal of Algebraic Geometry}, 23(1):63--90, 2014.

\bibitem{burgisser2013algebraic}
P.~B{\"u}rgisser, M.~Clausen, and M.~A. Shokrollahi.
\newblock {\em Algebraic complexity theory}, volume 315.
\newblock Springer Science \& Business Media, 2013.

\bibitem{chen2019numerical}
J.~Chen and J.~Kileel.
\newblock Numerical implicitization.
\newblock {\em Journal of Software for Algebra and Geometry}, 9(1):55--63,
  2019.

\bibitem{christandl2019barriers}
M.~Christandl, P.~Vrana, and J.~Zuiddam.
\newblock Barriers for fast matrix multiplication from irreversibility.
\newblock In {\em 34th Computational Complexity Conference (CCC 2019)}, 2019.

\bibitem{conner2022rank}
A.~Conner, F.~Gesmundo, J.~M. Landsberg, and E.~Ventura.
\newblock Rank and border rank of {K}ronecker powers of tensors and
  {S}trassen's laser method.
\newblock {\em computational complexity}, 31(1):1, 2022.

\bibitem{conner2021towards}
A.~Conner, F.~Gesmundo, J.~M. Landsberg, E.~Ventura, and Y.~Wang.
\newblock Towards a geometric approach to {S}trassen’s asymptotic rank
  conjecture.
\newblock {\em Collectanea mathematica}, 72(1):63--86, 2021.

\bibitem{duff2019solving}
T.~Duff, C.~Hill, A.~Jensen, K.~Lee, A.~Leykin, and J.~Sommars.
\newblock Solving polynomial systems via homotopy continuation and monodromy.
\newblock {\em IMA Journal of Numerical Analysis}, 39(3):1421--1446, 2019.

\bibitem{galkazka2023distinguishing}
M.~Ga{\l}{\k{a}}zka, T.~Ma{\'n}dziuk, and F.~Rupniewski.
\newblock Distinguishing secant from cactus varieties.
\newblock {\em Foundations of Computational Mathematics}, 23(4):1167--1214,
  2023.

\bibitem{griffin2014numerical}
Z.~A. Griffin, J.~D. Hauenstein, C.~Peterson, and A.~J. Sommese.
\newblock Numerical computation of the {H}ilbert function and regularity of a
  zero dimensional scheme.
\newblock In {\em Connections between algebra, combinatorics, and geometry},
  pages 235--250. Springer, 2014.

\bibitem{hauenstein2013equations}
J.~D. Hauenstein, C.~Ikenmeyer, and J.~M. Landsberg.
\newblock Equations for lower bounds on border rank.
\newblock {\em Experimental Mathematics}, 22(4):372--383, 2013.

\bibitem{hauenstein2010witness}
J.~D. Hauenstein and A.~J. Sommese.
\newblock Witness sets of projections.
\newblock {\em Applied Mathematics and Computation}, 217(7):3349--3354, 2010.

\bibitem{hauenstein2013membership}
J.~D. Hauenstein and A.~J. Sommese.
\newblock Membership tests for images of algebraic sets by linear projections.
\newblock {\em Applied Mathematics and Computation}, 219(12):6809--6818, 2013.

\bibitem{kaski2025universal}
P.~Kaski and M.~Micha{\l}ek.
\newblock A universal sequence of tensors for the asymptotic rank conjecture.
\newblock In {\em 16th Innovations in Theoretical Computer Science Conference
  (ITCS 2025)}, pages 64--1. Schloss Dagstuhl--Leibniz-Zentrum f{\"u}r
  Informatik, 2025.

\bibitem{michalek2021invitation}
M.~Micha{\l}ek and B.~Sturmfels.
\newblock {\em Invitation to nonlinear algebra}, volume 211.
\newblock American Mathematical Soc., 2021.

\bibitem{schonhage1981partial}
A.~Sch{\"o}nhage.
\newblock Partial and total matrix multiplication.
\newblock {\em SIAM Journal on Computing}, 10(3):434--455, 1981.

\bibitem{sommese2001using}
A.~J. Sommese, J.~Verschelde, and C.~W. Wampler.
\newblock Using monodromy to decompose solution sets of polynomial systems into
  irreducible components.
\newblock In {\em Applications of algebraic geometry to coding theory, physics
  and computation}, pages 297--315. Springer, 2001.

\bibitem{sommese2005numerical}
A.~J. Sommese, C.~W. Wampler, et~al.
\newblock {\em The numerical solution of systems of polynomials arising in
  engineering and science}.
\newblock World Scientific, 2005.

\bibitem{strassen1986asymptotic}
V.~Strassen.
\newblock The asymptotic spectrum of tensors and the exponent of matrix
  multiplication.
\newblock In {\em 27th Annual Symposium on Foundations of Computer Science
  (sfcs 1986)}, pages 49--54. IEEE, 1986.

\bibitem{strassen1988asymptotic}
V.~Strassen.
\newblock The asymptotic spectrum of tensors.
\newblock {\em Journal f{\"u}r die reine und angewandte Mathematik},
  390:102--152, 1988.

\bibitem{wigderson2022asymptotic}
A.~Wigderson and J.~Zuiddam.
\newblock Asymptotic spectra: Theory, applications and extensions.
\newblock {\em Manuscript}, 2022.

\end{thebibliography}
\end{document}